\documentclass[11pt]{article}
\usepackage{enumerate}
\usepackage{amssymb,a4wide,latexsym,makeidx,epsfig,fleqn}
\usepackage{amsthm}
\usepackage{amsmath}
\usepackage{enumerate}
\newpage
\newtheorem{theorem}{Theorem}[section]

\newtheorem{definition}[theorem]{Definition}
\newtheorem{lemma}[theorem]{Lemma}

\def\Tr{\mathop{\rm Tr }\nolimits}
\def\tr{\mathop{\rm tr }\nolimits}
\def\Res{\mathop{\rm Res }\nolimits}

\begin{document}
\textwidth 150mm \textheight 225mm
\title{The characteristic polynomials of uniform hypercycles with length four
\thanks{Supported by the National Natural Science Foundation of China (Nos. 12271439 and 12301452).}}
\author{{Cunxiang Duan$^{a}$\footnote{Corresponding author.}, Ligong Wang$^{b}$, Yulong Wei$^{a}$ }\\
{\small $^{a}$School of Mathematics, Taiyuan University of Technology, Taiyuan, Shanxi 030024, P.R. China}\\
{\small $^{b}$School of Mathematics and Statistics, Northwestern
Polytechnical University, Xi'an 710129, P.R. China}\\
\\{\small E-mail: cxduanmath@163.com; lgwangmath@163.com;weiyulong@tyut.edu.cn}\\}
\date{}
\maketitle
\begin{center}
\begin{minipage}{120mm}
\vskip 0.3cm
\begin{center}
{\small {\bf Abstract}}
\end{center}
{\small Let $C_{m}$ be a cycle with length $m.$ The $k$-uniform hypercycle with length $m$ obtained by adding $k-2$ new vertices in every edge of $C_{m},$ denoted by $C_{m,k}.$ In this paper, we obtain some trace formulas of uniform hypercycles with length four. Moreover, we give the characteristic polynomials of uniform hypercycles with length four.
}
\end{minipage}
\end{center}

\section{Introduction }
\label{sec:ch6-introduction}

Let $\mathcal{T}=(t_{i_{1}i_{2}\ldots i_{k}})$ be a $k$-order and $n$-dimensional tensor over copmplex $\mathbb{C},$ that is, a multidimensional array, $ 1\leq i_{1},i_{2},\ldots,i_{k}\leq n.$ For a vertex set $V,$ if any edge $e$ in the edge set $E$ is a subset of $V,$ then $H=(V,E)$ is called a hypergraph. Further, $H$ is $k$-uniform if every edge $e\in E$ such that $|e|=k.$ Note that the $2$-uniform hypergraph is a graph. Recently, the study of spectral hypergraph via eigenvalues of tensors has attracted the attention and research of many researchers. At the same time, it has also achieved rich results about the spectral theory of hypergraphs \cite{Bret,DW,GCH,KLQY}.

The characteristic polynomials of hypergraphs are an important research topic in spectra of hypergraphs. At present, one of the most commonly used method to study the characteristic polynomials of hypergraphs is the Poisson Formula. In 2015, Cooper and Dutle \cite{CoD}, by using the Poisson Formula, gave the spectra of "all ones" tensors and the characteristic polynomials of adjacency tensors of 3-uniform hyperstars. In 2018, Bao et al. \cite{BFWZ}, by using the Poisson Formula, obtained the characteristic polynomials of adjacency tensors of $k$-uniform hyperstars. In 2019, Chen and Bu \cite{CB}, by using the Poisson Formula, obtained a reduced formula of the characteristic polynomials of adjacency tensors of $k$-uniform hypergraphs with pendant edges. Further, they presented the characteristic polynomials of adjacency tensors of $k$-uniform hyperpaths. In 2021, Zheng \cite{Z} gave the characteristic polynomials of adjacency tensors of complete 3-uniform hypergraphs.

Moreover, some scholars have also studied the characteristic polynomial coefficients of uniform hypergraphs. In 2012, Cooper and Dutle \cite{CoDu} researched the characteristic polynomial first $k+1$ coefficients of the adjacency tensor of the $k$-uniform hypergraph and the characteristic polynomial of the adjacency tensor of a single hyperedge by using the characteristic equation, respectively. In 2014, Zhou et al. \cite{ZSWB} gave the first $k$ coefficient expression of the characteristic polynomial of the signless (Laplacian) tensor of a $k$-uniform hypergraph. Note that it is not easy to get an explicit characteristic polynomial expression of uniform hypergraphs even if some special uniform hypergraphs.

In fact, the spectra of hypergraphs can be studied by tensor traces. In 2015, Shao, Qi and Hu \cite{SQH} gave some new trace formulas and obtained the characterization of $k$-symmetric spectra of adjacency tensors of $k$-uniform hypergraphs. In 2021, Clark and Cooper \cite{CC} presented the Harary-Sachs theorem of $k$-uniform hypergraphs, which generalized the result of graphs. More results, see \cite{CBZ,HHLQ}. Motivated by above papers, we mainly consider the characteristic polynomials of adjacency tensors of uniform hypercycles with length four.

Let $H=(V,E)$ be a $k$-uniform hypergraph, where $V=\{v_1, v_{2}, \ldots, v_{m(k-1)+1}\}$ and $E=\{e_1, e_2, \ldots,e_m\}.$ If $e_i=\{v_{(i-1)(k-1)+1},v_{(i-1)(k-1)+2},\ldots, v_{i(k-1)+1}\}\in E$ for $i=1,2,\ldots,m,$  and  $v_{1}=v_{m(k-1)+1},$ then $H$ is a $k$-uniform hypercycle with length $m,$ denoted by $C_{m,k}.$

In Section 2, we introduce some basic definitions and properties of tensors and hypergraphs. We also give some useful lemmas which will be used in Section 3. In Section 3, we give some traces of adjacency tensors of $k$-uniform hypercycles with length four, and present the characteristic polynomials of adjacency tensors of $k$-uniform hypercycles with length four by using these traces.

\section{Preliminaries}
\label{sec:ch-sufficient}
In this section, we mainly introduce some basic definitions and properties of tensors and hypergraphs. Further, we present some useful Lemmas.

For a $k$-order and $n$-dimensional tensor $\mathcal{T}=(t_{i_{1}i_{2}\ldots i_{k}})$ and a vector $x =(x_{1}, x_{2}, \ldots, x_{n})^{\top}$ over $\mathbb{C}^{n}$, $\mathcal{T}x$ is a vector and
$$(\mathcal{T}x)_{i}=\sum\limits^{n}_{i_{2},\ldots,i_{k}=1}t_{ii_{2}\ldots i_{k}}x_{i_{2}}\cdots x_{i_{k}},~1\leq i \leq n.$$

\noindent\begin{definition}\label{de:c1} \cite{Qi}
Let $\mathcal{T}$ be a $k$-order and $n$-dimensional nonzero tensor, $x\in \mathbb{C}^{n}$ be a nonzero vector and $x^{[k-1]}=(x_{1}^{k-1}, x_{2}^{k-1}, \ldots, x_{n}^{k-1})^{\top}$. If there exists a number $\lambda\in \mathbb{C}$ and $$\mathcal{T} x=\lambda x^{[k-1]},$$ then $\lambda$ and $x$ respectively is an eigenvalue of $\mathcal{T}$ and an eigenvector of $\mathcal{T}$ corresponding to $\lambda.$
\end{definition}
For a $k$-order and $n$-dimensional tensor $\mathcal{T}=(t_{i_{1}i_{2}\cdots i_{k}}),$ the characteristic polynomial $\phi_{\mathcal{T}}(\lambda)$ of $\mathcal{T}$ is the resultant $\Res(\lambda x^{[k-1]}-\mathcal{T} x^{k-1}),$ and $\phi_{\mathcal{T}}(\lambda)$ is a monic polynomial in $\lambda$ of degree $n(k-1)^{n-1}.$ The $j$-th order trace $\Tr_{j}(\mathcal{T})$ \cite{MS} of $\mathcal{T}$ is $$\Tr_{j}(\mathcal{T})=(k-1)^{n-1}\sum_{j_{1}+j_{2}+\cdots+j_{n}=j}[\prod_{i=1}^{n}\frac{1}{(j_{i}(k-1))!}\big(\sum_{l \in [n]^{k-1}}t_{il}\frac{\partial}{\partial x_{il}}\big)^{j_{i}}] \tr (X^{j(k-1)}),$$
where the auxiliary $n\times n$ matrix $X = (x_{ij})$, $\frac{\partial}{\partial x_{il}}$ equals $ \frac{\partial}{\partial x_{il_{2}}}\frac{\partial}{\partial x_{il_{3}}}\cdots \frac{\partial}{\partial x_{il_{k}}}$ for $l = l_{2}\cdots l_{k},$ and  $j_{1}, \ldots, j_{n}$ run over all nonnegative integers with $j_{1} +\cdots+j_{n} =j.$ Moreover, Morozov and Shakirov \cite{MS} gave a formula for calculating the characteristic polynomial of $\mathcal{T}$ by using "Schur polynomials" in the generalized traces of $\mathcal{T},$ that is,
$$\phi_{\mathcal{T}}(\lambda)=\sum_{j=0}^{s}P_{j}(-\frac{\Tr_{1}(\mathcal{A})}{1}, -\frac{\Tr_{2}(\mathcal{A})}{2}, \ldots, -\frac{\Tr_{j}(\mathcal{A})}{j})\lambda^{s-j},$$
where $s=n(k-1)^{n-1}$ and the Schur function $P_{j}(p_{1},  \ldots, p_{j})=\sum\limits_{i=1}^{j}\sum\limits_{h_{1}+h_{2}+\cdots+h_{i}=j}\frac{p_{h_{1}} \cdots p_{h_{i}}}{i!}~(P_{0}=1).$ Note that the $j$-th order trace of $\mathcal{T}$ is the sum of $j$ power of all eigenvalues of $\mathcal{T}$ \cite{HHLQ}.

\noindent\begin{definition}\label{de:c1} \cite{CoDu}
Let $H = (V, E)$ be a $k$-uniform hypergraph with $n$ vertices. The adjacency tensor of $H$ is an $n$-dimensional tensor $\mathcal{A}_{H}=(a_{i_{1}i_{2}\ldots i_{k}})$ of order $k$ and $$a_{i_{1}i_{2}\ldots i_{k}}=\left\{
\begin{array}{ll}
\frac{1}{(k-1)!},& \mbox {if}   ~\{i_{1},i_{2},\ldots, i_{k}\} \in E,
\\
0,& \mbox {otherwise}.
\end{array}
\right.$$
\end{definition}
The degree $d_{v}$ of a $k$-uniform hypergraph $H$ equals $|e_{v}|$, where $e_{v}$ is a set of edges that incident with the vertex $v.$ If $d_{v}$ is the multiple of $k,$ for all vertex $v\in V,$ then $H$ is called $k$-valent. Let $[n]=\{1,2,\ldots,n\}.$ For a positive integer $j,$ we define $$\mathcal{F}_{j}=\{(i_{1}\alpha_{1},\ldots,i_{j}\alpha_{j}) \mid 1 \leq i_{1} \leq i_{2} \leq\cdots \leq i_{j}\leq n,\alpha_{1}, \ldots, \alpha_{j}\in [n]^{k-1}\}.$$ Let denote $\pi_{F}(\mathcal{T})=t_{i_{1}\alpha_{1}}\cdots t_{i_{j}\alpha_{j}}$ for $F=(i_{1}\alpha_{1},\ldots,i_{j}\alpha_{j})\in \mathcal{F}_{j}$ and a tensor $\mathcal{T}=(t_{i_{1}i_{2}\cdots i_{k}}).$ For a digraph $D=(V,A)$, the in-degree (out-degree) of $v_{i}$ in $D$ is the number of arcs incident to (from) $v_{i},$ denoted by $d^{-}_{v_{i}}~(d^{+}_{v_{i}}).$
\noindent\begin{definition}\label{de:4-1} \cite{SQH}
For $F=(i_{1}\alpha_{1},\ldots,i_{j}\alpha_{j})\in \mathcal{F}_{j},$ let $E(F)=\bigcup\limits_{h=1}^{j}E_{h}(F),$ where $E_{h}(F)=\{(i_{h},v_{1}), \ldots, (i_{h},v_{k-1}) \}$ denotes the arc multi-set if $\alpha_{h}=v_{1}\ldots v_{k-1}.$ Let $D=(V,A)$ denote the (multi-)digraph corresponding to $E(F),$ where $V$ and $A$ is the vertex set and arc set, respectively. Then
\\$(1)$ $b(F)=\prod\limits_{a\in A} m(a)!$ and $c(F)=\prod\limits_{v\in V}d^{+}_{v}!,$ where $m(a)$ is the multiplicity of the arc $a.$
\\$(2)$ Denote $W(F)$ the set of all Eulerian closed walks with $E(F).$
\end{definition}
Note that multiple arcs of $W(F)$ are not distinguished and $W(F)=\emptyset$ if $F$ is not $k$-valent.

\noindent\begin{lemma}\label{le:4-4}\cite{SQH}
Let $\mathcal{T}=(t_{i_{1}i_{2}\cdots i_{k}})$ be an $n$-dimensional tensor of order $k$. Then
$$\Tr_{j}(\mathcal{T})=(k-1)^{n-1}\sum_{F\in \mathcal{F}'_{j}}\frac{b(F)}{c(F)}\pi_{F}(\mathcal{T})\lvert W(F)\rvert,$$
where $\mathcal{F}'_{j}=\{F\in \mathcal{F}_{j} \mid F$ is $k$-valent$\}.$
\end{lemma}

Let $G=(V,E)$ be a graph. If each edge of $G$ adds $k-2$ new vertices, then the $k$-unifrom hypergraph is $k$-power of $G,$ denoted by $G^{k}.$ If the edge sign function $\pi: E\rightarrow \{+1,-1\},$ then $(G,\pi)$ is a signed graph, denoted by $G_{\pi}.$ If a graph is an (induced) subgraph of $G_{\pi},$ then the garaph is a signed (induced) subgraphs of $G.$
\noindent\begin{lemma}\label{le:4-5}\cite{CDB} $\lambda\in \mathbb{C}$ is an eigenvalue of $G^{k}$ if and only if
\\$(1)$ for $k = 3,$ $\beta$ is an eigenvalue of some signed induced subgraph of $G$ and $\beta^{2}= \lambda^{k}$;
\\$(2)$ for $k\geq 4,$ $\beta$ is an eigenvalue of some signed subgraph of $G$ and $\beta^{2}= \lambda^{k}.$
\end{lemma}

\section{The characteristic polynomial of $k$-uniform hypercycles with length four}
\label{sec:ch-inco}
In this section, we present the characteristic polynomials of $k$-uniform hypercycles with length four. Note that the number of spanning trees $\tau(D)$ and the number of Eulerian cycles $\varepsilon(D)$ of a digraph $D$ can be calculate by using the Matrix-Tree Theorem \cite{CvRS} and BEST Theorem \cite{AB}, respectively.

\noindent\begin{theorem}\label{th:3-1}
Let $C_{4,k}$ be a $k(\geq 3)$-uniform hypercycle with length four, and $\mathcal{A}_{C_{4,k}}$ be its adjacency tensor. Then $$\Tr_{k}(\mathcal{A}_{C_{4,k}})=4k^{k-1}(k-1)^{3k-4}.$$
\end{theorem}
\noindent\textbf{Proof.} By Lemma \ref{le:4-4}, we consider $F=(i_{1}\alpha_{1}, \ldots, i_{k}\alpha_{k})\in \mathcal{F}'_{k}.$ If $\pi_{F}(\mathcal{A}_{C_{4,k}})\neq 0,$ then we know that all elements of $F$ correspond to edges of $C_{4,k}.$ Since $F\in \mathcal{F}'_{k}$ is $k$-valent, each vertex of the edge in $F$ occurs the times that is the multiple of $k.$ Thus, all elements $i_{h}\alpha_{h}$ of $F$ only correspond to some edge of $C_{4,k},$ $1 \leq h \leq k.$ If $|W(F)|\neq 0,$ then the out-degree is equal to in-degree of every vertex in the directed graph corresponding to $F$, that is, every vertex as the first entry occurs one time in $F$. For an edge of $C_{4,k},$ the total number of such $F$ is $[(k-1)!]^{k}.$ For each such $F,$ $E(F)$ induces a complete digraph $D_{1}$ on $k$ vertices. Hence, $$b(F)=1, c(F)= [(k-1)!]^{k} , ~\pi_{F}(\mathcal{A}_{C_{4,k}})=\big[\frac{1}{(k-1)!}\big]^{k}.$$ By the Matrix-Tree Theorem \cite{CvRS}, we have $\tau(D_{1})=k^{k-2}$. By BEST Theorem \cite{AB}, we have $\varepsilon(D_{1})=[(k-2)!]^{k}k^{k-2}$. Because each Eulerian cycle has $k(k-1)$ arcs, we have $$\lvert W(F)\rvert=k(k-1)[(k-2)!]^{k}k^{k-2}.$$ Since $C_{4,k}$ has four edges, we have
\begin{align*}
\Tr_{k}(\mathcal{A}_{C_{4,k}})&=4(k-1)^{n-1}\cfrac{[(k-1)!]^{k}}{[(k-1)!]^{k}}\big[\cfrac{1}{(k-1)!}\big]^{k}k(k-1)[(k-2)!]^{k}k^{k-2}
\\&=4k^{k-1}(k-1)^{n-k}=4k^{k-1}(k-1)^{3k-4}.
\end{align*}

\noindent\begin{theorem}\label{th:3-2}
Let $C_{4,k}$ be a $k(\geq 3)$-uniform hypercycle with length four, and $\mathcal{A}_{C_{4,k}}$ be its adjacency tensor. Then $$\Tr_{2k}(\mathcal{A}_{C_{4,k}})=4k^{k-1}(k-1)^{3k-4}+8k^{2k-3}(k-1)^{2k-3}.$$
\end{theorem}
\noindent\textbf{Proof.}
By Lemma \ref{le:4-4}, we consider $F=(i_{1}\alpha_{1}, \ldots, i_{2k}\alpha_{2k})\in \mathcal{F}'_{2k}.$ If $\pi_{F}(\mathcal{A}_{C_{4,k}})\neq 0,$ then we know that all elements of $F$ correspond to edges of $C_{4,k}.$ since $F\in \mathcal{F}'_{2k}$ is $k$-valent, each vertex of the edge in $F$ occurs the times that is the multiple of $k.$ Thus, $F$ has the following two cases.

{\bf Case 1.} All elements $i_{h}\alpha_{h}$ of $F$ correspond to some edge of $C_{4,k},$ $1 \leq h \leq 2k.$

Similar to Theorem \ref{th:3-1}, we know that $E(F)$ induces a complete multi-digraph $D'_{1}$ on $k$ vertices and the multiplicity of each arc of $E(F)$ is 2. Hence, the total number of such $F$ is $[(k-1)!]^{2k},$  and $$b(F)=(2!)^{k(k-1)}, c(F)= [(2(k-1))!]^{k} , ~\pi_{F}(\mathcal{A}_{C_{4,k}})=\big[\frac{1}{(k-1)!}\big]^{2k}.$$ By the Matrix-Tree Theorem \cite{CvRS}, we have $\tau(D'_{1})=2^{k-1}k^{k-2}$. By BEST Theorem \cite{AB}, we have $\varepsilon(D'_{1})=[(2(k-1)-1)!]^{k}2^{k-1}k^{k-2}$. Because each Eulerian cycle has $2k(k-1)$ arcs and multi-arcs of $W(F)$ are not labelled, we have $$\lvert W(F)\rvert=\cfrac{2k(k-1)[(2(k-1)-1)!]^{k}2^{k-1}k^{k-2}}{(2!)^{k(k-1)}}.$$ Consider all edges, for all such $F,$ we know that the total contribution to $\Tr_{2k}(\mathcal{A}_{C_{4,k}})$ is
\begin{align*}
&4(k-1)^{n-1}\cfrac{[(k-1)!]^{2k}(2!)^{k(k-1)}}{[(2(k-1))!]^{k}}\big[\cfrac{1}{(k-1)!}\big]^{2k}\cfrac{2k(k-1)[(2(k-1)-1)!]^{k}2^{k-1}k^{k-2}}{(2!)^{k(k-1)}}
\\&=4k^{k-1}(k-1)^{n-k}=4k^{k-1}(k-1)^{3k-4}.
\end{align*}

{\bf Case 2.} All elements $i_{h}\alpha_{h}$ of $F$ correspond to two edges of $C_{4,k},$ $1 \leq h \leq 2k.$

Since $F\in \mathcal{F}'_{2k}$ is $k$-valent, we know that $k$ elements of $F$ correspond to the same edge and $k$ elements of $F$ correspond to the same other edge. If $|W(F)|\neq 0,$ then every vertex of each edge as the first entry occurs the same number of times in $F$. We only consider two incident edges since the multi-digraph corresponding to $E(F)$ is connected. For two incident edges of $C_{4,k},$ the number of orderings for the first entry is 2, and the number of orderings for $\alpha_{h}$ is $[(k-1)!]^{2k}$. The total number of such $F$ is $2[(k-1)!]^{2k}.$ Thus, $E(F)$ induces a digraph $D_{2}$ on $2k-1$ vertices. Since $D_{2}$ must be connected, we know $F=(i_{1}\alpha_{1}, \ldots,i_{k}\alpha_{k},i_{k}\alpha_{k+1},\ldots,i_{2k-1}\alpha_{2k-1}).$ Hence, $$b(F)=1,~c(F)= [(k-1)!]^{2k-2}(2(k-1))!,~ \pi_{F}(\mathcal{A}_{C_{4,k}})=\big[\cfrac{1}{(k-1)!}\big]^{2k},$$ By the Matrix-Tree Theorem \cite{CvRS}, we have $\tau(D_{2})=k^{2k-4}$. By BEST Theorem \cite{AB}, we have $\varepsilon(D_{2})=[(k-2)!]^{2k-2}(2(k-1)-1)!k^{2k-4}$. Because each Eulerian cycle has $2k(k-1)$ arcs, we have $$\lvert W(F)\rvert=2k(k-1)[(k-2)!]^{2k-2}(2(k-1)-1)!k^{2k-4}.$$ Consider all edges, for all such $F,$ we know that the total contribution to $\Tr_{2k}(\mathcal{A}_{C_{4,k}})$ is
\begin{align*}
&4(k-1)^{n-1}\cfrac{2[(k-1)!]^{2k}}{[(k-1)!]^{2k-2}(2(k-1))!}\big[\cfrac{1}{(k-1)!}\big]^{2k} 2k(k-1)[(k-2)!]^{2k-2}(2(k-1)-1)!k^{2k-4}
\\&=8k^{2k-3}(k-1)^{n-2k+1}=8k^{2k-3}(k-1)^{2k-3}.
\end{align*}

Thus, we have
\begin{align*}
\Tr_{2k}(\mathcal{A}_{C_{4,k}})=4k^{k-1}(k-1)^{3k-4}+8k^{2k-3}(k-1)^{2k-3}.
\end{align*}

\noindent\begin{theorem}\label{th:3-3}
Let $C_{4,k}$ be a $k(\geq 3)$-uniform hypercycle with length four, and $\mathcal{A}_{C_{4,k}}$ be its adjacency tensor. Then $$\Tr_{3k}(\mathcal{A}_{C_{4,k}})=4k^{k-1}(k-1)^{3k-4}+24k^{2k-3}(k-1)^{2k-3}+12k^{3k-5}(k-1)^{k-2}.$$
\end{theorem}
\noindent\textbf{Proof.} By Lemma \ref{le:4-4}, we consider $F=(i_{1}\alpha_{1}, \ldots, i_{3k}\alpha_{3k})\in \mathcal{F}'_{3k}.$ If $\pi_{F}(\mathcal{A}_{C_{4,k}})\neq 0,$ then we know that all elements of $F$ correspond to edges of $C_{4,k}.$ Since $F\in \mathcal{F}'_{3k}$ is $k$-valent, each vertex of the edge in $F$ occurs the times that is the multiple of $k.$  Thus, $F$ has the following three cases.

{\bf Case 1.} All elements $i_{h}\alpha_{h}$ of $F$ correspond to some edge of $C_{4,k},$ $1 \leq h \leq 3k.$

Similar to Case 1 of Theorem \ref{th:3-2}, we know that $E(F)$ corresponds a complete multi-digraph with $k$ vertices and the multiplicity of each arc of $E(F)$ is 3. Hence, we know that the total number of such $F$ is $[(k-1)!]^{3k},$ and $$b(F)=(3!)^{k(k-1)}, c(F)= [(3(k-1))!]^{k} , ~\pi_{F}(\mathcal{A}_{C_{4,k}})=\big[\frac{1}{(k-1)!}\big]^{3k},$$ $$\lvert W(F)\rvert=\cfrac{3k(k-1)[(3(k-1)-1)!]^{k}3^{k-1}k^{k-2}}{(3!)^{k(k-1)}}.$$ Since $C_{4,k}$ has four edges, for all such $F,$ we know that the total contribution to $\Tr_{3k}(\mathcal{A}_{C_{4,k}})$ is
\begin{align*}
&4(k-1)^{n-1}\cfrac{[(k-1)!]^{3k}(3!)^{k(k-1)}}{[(3(k-1))!]^{k}}\big[\cfrac{1}{(k-1)!}\big]^{3k}\cfrac{3k(k-1)[(3(k-1)-1)!]^{k}3^{k-1}k^{k-2}}{(3!)^{k(k-1)}}
\\&=4k^{k-1}(k-1)^{n-k}=4k^{k-1}(k-1)^{3k-4}.
\end{align*}

{\bf Case 2.} All elements $i_{h}\alpha_{h}$ of $F$ correspond to two edges of $C_{4,k},$ $1 \leq h \leq 3k.$

Assume that $ak$ elements of $F$ correspond to the same edge and $bk$ elements of $F$ correspond to the same other edge, $a+b=3.$ If $|W(F)|\neq 0,$ then the in-degree is equal to out-degree of every vertex of the multi-digraph corresponding to $F.$ We only consider two incident edges since the multi-digraph corresponding to $E(F)$ is connected. For two incident edges of $C_{4,k},$ the number of orderings for the first entry is $\binom{t}{a}$, and the number of orderings for $\alpha_{h}$ is $[(k-1)!]^{3k}$. The total number of such $F$ is $\sum\limits_{a=1}^{2}\binom{3}{a}[(k-1)!]^{3k}.$

We know that $E(F)$ induces a multi-digraph $D'_{2}$ on $2k-1$ vertices. Since $D'_{2}$ must be connected, $F$ is an appropriate ordering of $(i_{1}\alpha^{1}_{1}, \ldots,i_{1}\alpha^{a}_{1},\ldots,i_{k}\alpha^{1}_{k},\ldots,i_{k}\alpha^{a}_{k},i_{k}\alpha^{1'}_{k}\ldots,i_{k}\alpha^{b'}_{k},\cdots,\\ i_{2k-1}\alpha^{1'}_{2k},\cdots,i_{2k-1}\alpha^{b'}_{2k+1},\ldots, i_{3k-2}\alpha^{b'}_{3k}),$ where $\alpha^{j}_{h}$ (resp. $\alpha^{l'}_{h}$) has the same elements regardless of ordering, for $j=1,2,\cdots,a,$ $l=1,2, \cdots, b.$ Hence, $$b(F)=(a!)^{k(k-1)}(b!)^{k(k-1)},~c(F)= [(a(k-1))!]^{k-1}[(b(k-1))!]^{k-1}(3(k-1))!,$$ and $\pi_{F}(\mathcal{A}_{C_{4,k}})=\big[\cfrac{1}{(k-1)!}\big]^{3k}.$ By the Matrix-Tree Theorem \cite{CvRS}, we have $\tau(D'_{2})=a^{k-1}b^{k-1}k^{2k-4}$. By BEST Theorem \cite{AB}, we have $\varepsilon(D'_{2})=[(a(k-1)-1)!]^{k}[(b(k-1)-1)!]^{k}((a+b)(k-1)-1)!a^{k-1}b^{k-1}k^{2k-4}.$ Because each Eulerian cycle has $3k(k-1)$ arcs and multi arcs of $W(F)$ are not labelled, we have $$\lvert W(F)\rvert=\cfrac{3k(k-1)[(a(k-1)-1)!]^{k-1}[(b(k-1)-1)!]^{k-1}(3(k-1)-1)!a^{k-1}b^{k-1}k^{2k-4}}{(a!)^{k(k-1)}(b!)^{k(k-1)}}.$$  Consider all edges, for all such $F,$ we know that the total contribution to $\Tr_{3k}(\mathcal{A}_{C_{4,k}})$ is
\begin{align*}
&4(k-1)^{n-1}\cfrac{\sum_{a=1}^{2}\binom{3}{a}[(k-1)!]^{3k}(a!)^{k(k-1)}(b!)^{k(k-1)}}{[(a(k-1))!]^{k-1}[(b(k-1))!]^{k-1}(3(k-1))!}\big[\cfrac{1}{(k-1)!}\big]^{3k} \\&~~~~\cfrac{3k(k-1)[(a(k-1)-1)!]^{k-1}[(b(k-1)-1)!]^{k-1}(3(k-1)-1)!a^{k-1}b^{k-1}k^{2k-4}}{(a!)^{k(k-1)}(b!)^{k(k-1)}}
\\&=24k^{2k-3}(k-1)^{n-2k+1}=24k^{2k-3}(k-1)^{2k-3}.
\end{align*}

{\bf Case 3.} All elements $i_{h}\alpha_{h}$ of $F$ correspond to three edges of $C_{4,k},$ $1 \leq h \leq 3k.$

If $|W(F)|\neq 0$, then the in-degree is equal to out-degree of every vertex of the multi-digraph corresponding to $F.$ Thus, $F$ is an appropriate ordering of $(i_{1}\alpha_{1}, i_{2}\alpha_{2},\ldots,i_{k}\alpha_{k},i_{k}\alpha_{k+1},\ldots,\\i_{2k-1}\alpha_{2k},i_{2k-1}\alpha_{2k+1},\ldots, i_{3k-2}\alpha_{3k}).$ We know that $E(F)$ induces a digraph $D_{3}$ on $3k-2$ vertices. For each such $F,$ the number of orderings for the first entry is 4, and the number of orderings for the $\alpha_{h}$ is $[(k-1)!]^{3k}$. Hence, the total number of such $F$ is $4[(k-1)!]^{3k},$ and $$b(F)=1, ~c(F)= [(k-1)!]^{3k-4}[(2(k-1))!]^{2},~\pi_{F}(\mathcal{A}_{C_{4,k}})=\big[\cfrac{1}{(k-1)!}\big]^{3k}.$$ By the Matrix-Tree Theorem \cite{CvRS}, we have $\tau(D_{3})=k^{3k-6}.$ By BEST Theorem \cite{AB}, we have $\varepsilon(D_{3})=[(k-2)!]^{3k-4}[(2(k-1)-1)!]^{2}k^{3k-6}.$
Because each Eulerian cycle has $3k(k-1)$ arcs, we have
\begin{align*}
\lvert W(F)\rvert&=3k^{3k-5}(k-1)[(k-2)!]^{3k-4}[(2(k-1)-1)!]^{2}.
\end{align*}
Consider all edges, for all such $F,$ we know that the total contribution to $\Tr_{3k}(\mathcal{A}_{C_{4,k}})$ is
\begin{align*}
&4(k-1)^{n-1}\cfrac{4[(k-1)!]^{3k}}{[(k-1)!]^{3k-4}[(2(k-1))!]^{2}}\big[\cfrac{1}{(k-1)!}\big]^{3k}3k^{3k-5}(k-1)[(k-2)!]^{3k-4}[(2(k-1)-1)!]^{2}
\\&=12k^{3k-5}(k-1)^{n-3k+2}=12k^{3k-5}(k-1)^{k-2}.
\end{align*}

Thus, we have
\begin{align*}
\Tr_{3k}(\mathcal{A}_{C_{4,k}})=4k^{k-1}(k-1)^{3k-4}+24k^{2k-3}(k-1)^{2k-3}+12k^{3k-5}(k-1)^{k-2}
\end{align*}

\noindent\begin{theorem}\label{th:3-4}
Let $C_{4,k}$ be a $k(\geq 3)$-uniform hypercycle with length four, and $\mathcal{A}_{C_{4,k}}$ be its adjacency tensor. Then $$\Tr_{4k}(\mathcal{A}_{C_{4,k}})=4k^{k-1}(k-1)^{3k-4}+56k^{2k-3}(k-1)^{2k-3}+64k^{3k-5}(k-1)^{k-2}+40k^{4k-8}.$$
\end{theorem}

\noindent\textbf{Proof.} By Lemma \ref{le:4-4}, we consider $F=(i_{1}\alpha_{1}, \ldots, i_{4k}\alpha_{4k})\in \mathcal{F}'_{4k}.$ If $\pi_{F}(\mathcal{A}_{C_{4,k}})\neq 0,$ then all elements of $F$ correspond to edges of $C_{4,k}.$ since $F\in \mathcal{F}'_{4k}$ is $k$-valent, each vertex of the edge in $F$ occurs the times that is the multiple of $k.$ Thus, $F$ has the following four cases.

{\bf Case 1.} All elements $i_{h}\alpha_{h}$ of $F$ correspond to some edge of $C_{4,k},$ $1 \leq h \leq 4k.$

Similar to Case 1 of Theorem \ref{th:3-2}, we know that $E(F)$ corresponds a complete multi-digraph with $k$ vertices and the multiplicity of each arc of $E(F)$ is 4. Hence, the total number of such $F$ is $[(k-1)!]^{4k},$ and $$b(F)=(4!)^{k(k-1)}, c(F)= [(4(k-1))!]^{k} , ~\pi_{F}(\mathcal{A}_{C_{4,k}})=\big[\frac{1}{(k-1)!}\big]^{4k},$$ $$\lvert W(F)\rvert=\cfrac{4k(k-1)[(4(k-1)-1)!]^{k}4^{k-1}k^{k-2}}{(4!)^{k(k-1)}}.$$  Consider all edges, for all such $F,$ we know that the total contribution to $\Tr_{4k}(\mathcal{A}_{C_{4,k}})$ is
\begin{align*}
&4(k-1)^{n-1}\cfrac{[(k-1)!]^{4k}(4!)^{k(k-1)}}{[(4(k-1))!]^{k}}\big[\cfrac{1}{(k-1)!}\big]^{4k}\cfrac{4k(k-1)[(4(k-1)-1)!]^{k}4^{k-1}k^{k-2}}{(4!)^{k(k-1)}}
\\&=4k^{k-1}(k-1)^{n-k}=4k^{k-1}(k-1)^{3k-4}.
\end{align*}

{\bf Case 2.} All elements $i_{h}\alpha_{h}$ of $F$ correspond to two edges of $C_{4,k},$ $1 \leq h \leq 4k.$

Similar to Case 2 of Theorem \ref{th:3-3}, assume that $ak$ elements of $F$ correspond to the same edge and $bk$ elements of $F$ correspond to the same other edge, $a+b=4.$ We only consider two incident edges since the multi-digraph corresponding to $E(F)$ is connected. For two incident edges of $C_{4,k},$ the number of orderings for the first entry is $\binom{t}{a}$, and the number of orderings for the $\alpha_{h}$ is $[(k-1)!]^{4k}$. Thus, we know the total number of such $F$ is $\sum_{a=1}^{3}\binom{4}{a}[(k-1)!]^{4k},$ and $$b(F)=(a!)^{k(k-1)}(b!)^{k(k-1)}, \pi_{F}(\mathcal{A}_{C_{4,k}})=\big[\cfrac{1}{(k-1)!}\big]^{3k},$$ $$c(F)= [(a(k-1))!]^{k-1}[(b(k-1))!]^{k-1}(3(k-1))!,$$ $$\lvert W(F)\rvert=\cfrac{4k(k-1)[(a(k-1)-1)!]^{k-1}[(b(k-1)-1)!]^{k-1}(4(k-1)-1)!a^{k-1}b^{k-1}k^{2k-4}}{(a!)^{k(k-1)}(b!)^{k(k-1)}}.$$  Consider all edges, for all such $F,$ we know that the total contribution to $\Tr_{4k}(\mathcal{A}_{C_{4,k}})$ is
\begin{align*}
&4(k-1)^{n-1}\cfrac{\sum_{a=1}^{3}\binom{4}{a}[(k-1)!]^{4k}(a!)^{k(k-1)}(b!)^{k(k-1)}}{[(a(k-1))!]^{k-1}[(b(k-1))!]^{k-1}(4(k-1))!}\big[\cfrac{1}{(k-1)!}\big]^{4k} \\&~~~~\cfrac{4k(k-1)[(a(k-1)-1)!]^{k-1}[(b(k-1)-1)!]^{k-1}(4(k-1)-1)!a^{k-1}b^{k-1}k^{2k-4}}{(a!)^{k(k-1)}(b!)^{k(k-1)}}
\\&=56k^{2k-3}(k-1)^{n-2k+1}=56k^{2k-3}(k-1)^{2k-3}.
\end{align*}

{\bf Case 3.} All elements $i_{h}\alpha_{h}$ of $F$ correspond to three edges of $C_{4,k},$ $1 \leq h \leq 4k.$

Since $F\in \mathcal{F}'_{4k}$ is $k$-valent, there are $2k$ elements of $F$ corresponding to the same edge. If $|W(F)|\neq 0$, then the in-degree is equal to out-degree of every vertex of the multi-digraph corresponding to $F.$ Thus, $F$ is an appropriate ordering of $(i_{1}\alpha_{1}, i_{1}\alpha'_{1},\ldots,i_{k}\alpha_{k},i_{k}\alpha'_{k},i_{k}\alpha_{k+1},\\\ldots,i_{2k-1}\alpha_{2k-1},i_{2k-1}\alpha_{2k},\ldots, i_{3k-2}\alpha_{3k})$ or $(i_{1}\alpha_{1},\ldots,i_{k}\alpha_{k},i_{k}\alpha_{k+1},i_{k}\alpha'_{k+1},\ldots,i_{2k-1}\alpha_{2k},i_{2k-1}\alpha'_{2k},\\i_{2k-1}\alpha_{2k+1},\ldots, i_{3k-2}\alpha_{3k}),$ where $\alpha_{h}$ and $\alpha'_{h}$ have the same elements regardless of ordering, for $h=k+1,\ldots,2k.$

{\bf Subcase 3.1.} $F$ is an appropriate ordering of $(i_{1}\alpha_{1}, i_{1}\alpha'_{1},\ldots,i_{k}\alpha_{k},i_{k}\alpha'_{k},i_{k}\alpha_{k+1}\ldots,i_{2k-1}\alpha_{2k-1},\\i_{2k-1}\alpha_{2k},\ldots, i_{3k-2}\alpha_{3k}).$

For each such $F,$ the number of orderings for the first entry is 6, and the number of orderings for the $\alpha_{h}$ is $[(k-1)!]^{4k}$. Hence, the total number of such $F$ is $6[(k-1)!]^{4k},$ and $$b(F)=(2!)^{k}, ~c(F)= [(k-1)!]^{2k-3}[(2(k-1))!]^{k}(3(k-1))!,~\pi_{F}(\mathcal{A}_{C_{4,k}})=\big[\cfrac{1}{(k-1)!}\big]^{4k}.$$ Let $D_{3}$ be the multi-digraph induced by $E(F).$ By the Matrix-Tree Theorem \cite{CvRS}, we have $\tau(D_{3})=2^{k-1}k^{3k-6}$. By BEST Theorem \cite{AB}, we have $\varepsilon(D_{3})=2^{k-1}k^{3k-6}[(k-2)!]^{2k-3}[(2(k-1)-1)!]^{k}(3(k-1)-1)!.$
Because each Eulerian cycle has $4k(k-1)$ arcs and multi-arcs of $W(F)$ are not labelled, we have
\begin{align*}
\lvert W(F)\rvert=2k^{3k-5}(k-1)[(k-2)!]^{2k-3}[(2(k-1)-1)!]^{k}(3(k-1)-1)!.
\end{align*}
Consider all edges, for all such $F,$ we know that the total contribution to $\Tr_{4k}(\mathcal{A}_{C_{4,k}})$ is
\begin{align*}
&4(k-1)^{n-1}\cfrac{2\times6[(k-1)!]^{4k}(2!)^{k}}{[(k-1)!]^{2k-3}[(2(k-1))!]^{k}(3(k-1))!}\big[\cfrac{1}{(k-1)!}\big]^{4k}2k^{3k-5}(k-1)
\\&~~~~[(k-2)!]^{2k-3}[(2(k-1)-1)!]^{k}(3(k-1)-1)!
\\&=32k^{3k-5}(k-1)^{k-2}.
\end{align*}

{\bf Subcase 3.2.} $F$ is an appropriate ordering of $(i_{1}\alpha_{1},\ldots,i_{k}\alpha_{k},i_{k}\alpha_{k+1},i_{k}\alpha'_{k+1},\ldots,i_{2k-1}\alpha_{2k},\\i_{2k-1}\alpha'_{2k},i_{2k-1}\alpha_{2k+1},\ldots, i_{3k-2}\alpha_{3k}).$

For each such $F,$ the number of orderings for the first entry is 9, and the number of orderings for the $\alpha_{h}$ is $[(k-1)!]^{4k}$. Hence, the total number of such $F$ is $9[(k-1)!]^{4k},$ and $$b(F)=(2!)^{k},~c(F)= [(k-1)!]^{2k-2}[(2(k-1))!]^{k-2}[(3(k-1))!]^{2},~\pi_{F}(\mathcal{A}_{C_{4,k}})=\big[\cfrac{1}{(k-1)!}\big]^{4k}.$$
By calculating, we know $\lvert W(F)\rvert$ is the same as the above Subcase 3.1, i.e., $$\lvert W(F)\rvert=2k^{3k-5}(k-1)[(k-2)!]^{2k-2}[(2(k-1)-1)!]^{k-2}[(3(k-1)-1)!]^{2}.$$
Consider all edges, for all such $F,$ we know that the total contribution to $\Tr_{4k}(\mathcal{A}_{C_{4,k}})$ is
\begin{align*}
&4(k-1)^{n-1}\cfrac{9[(k-1)!]^{4k}(2!)^{k}}{[(k-1)!]^{2k-3}[(2(k-1))!]^{k}(3(k-1))!}\big[\cfrac{1}{(k-1)!}\big]^{4k}2k^{3k-5}(k-1)
\\&~~~~[(k-2)!]^{2k-2}[(2(k-1)-1)!]^{k-2}[(3(k-1)-1)!]^{2}
\\&=32k^{3k-5}(k-1)^{k-2}.
\end{align*}

{\bf Case 4.} All elements $i_{h}\alpha_{h}$ of $F$ correspond to all edges of $C_{4,k},$ $1 \leq h \leq 4k.$

If $|W(F)|\neq 0$, then the in-degree is equal to out-degree of every vertex of the (multi-)digraph corresponding to $F.$ Thus, $F$ is an appropriate ordering of $(i_{1}\alpha_{1}, i_{1}\alpha_{2},\ldots,i_{k}\alpha_{k+1},i_{k}\alpha_{k+2},\\\ldots,i_{2k-1}\alpha_{2k+1},i_{2k-1}\alpha_{2k+2},\ldots, i_{4k-4}\alpha_{4k})$ or $(i_{1}\alpha_{1},i_{1}\alpha'_{1},\ldots,i_{k}\alpha_{k},i_{k}\alpha'_{k},\ldots,\i_{2k-1}\alpha_{2k-1},\i_{2k-1}\alpha'_{2k-1},\\\ldots,i_{3k-2}\alpha_{3k-2},i_{3k-2}\alpha'_{3k-2},\ldots, i_{4k-4}\alpha_{4k-4}),$ where $\alpha_{h}$ and $\alpha'_{h}$ have the same elements regardless of ordering, for $h=1,k,2k-1,3k-2.$

{\bf Subcase 4.1.} $F$ is an appropriate ordering of $(i_{1}\alpha_{1}, i_{1}\alpha_{2},\ldots,i_{k}\alpha_{k+1},i_{k}\alpha_{k+2},\ldots,i_{2k-1}\alpha_{2k+1},\\i_{2k-1}\alpha_{2k+2},\ldots, i_{4k-4}\alpha_{4k}).$

For each such $F,$ the number of orderings for the first entry is 16, and the number of orderings for the $\alpha_{h}$ is $[(k-1)!]^{4k}$. Hence, the total number of such $F$ is $16[(k-1)!]^{4k},$ and $$b(F)=1, ~c(F)= [(k-1)!]^{4k-8}[(2(k-1))!]^{4},~\pi_{F}(\mathcal{A}_{C_{4,k}})=\big[\cfrac{1}{(k-1)!}\big]^{4k}.$$ Let $D_{4}$ be the multi-digraph induced by $E(F).$ By the Matrix-Tree Theorem \cite{CvRS}, we have $\tau(D_{4})=8k^{4k-9}.$ By BEST Theorem \cite{AB}, we have $\varepsilon(D_{4})=8k^{4k-9}[(k-2)!]^{4k-8}[(2(k-1)-1)!]^{4}.$
Because each Eulerian cycle has $4k(k-1)$ arcs, we have
\begin{align*}
\lvert W(F)\rvert&=32k^{4k-8}(k-1)[(k-2)!]^{4k-8}[(2(k-1)-1)!]^{4}.
\end{align*}
For all such $F,$ we know that the total contribution to $\Tr_{3k}(\mathcal{A}_{C_{4,k}})$ is
\begin{align*}
&(k-1)^{n-1}\cfrac{16[(k-1)!]^{4k}}{[(k-1)!]^{4k-8}[(2(k-1))!]^{4}}\big[\cfrac{1}{(k-1)!}\big]^{4k}32k^{4k-8}(k-1)
\\&[(k-2)!]^{4k-8}[(2(k-1)-1)!]^{4}
\\&=32k^{4k-8}.
\end{align*}

{\bf Subcase 4.2.} $F$ is an appropriate ordering of $(i_{1}\alpha_{1},i_{1}\alpha'_{1},\ldots,i_{k}\alpha_{k},i_{k}\alpha'_{k},\ldots,i_{2k-1}\alpha_{2k+1},\\i_{2k-1}\alpha'_{2k+1},\ldots,i_{3k-2}\alpha_{3k+2},i_{3k-2}\alpha'_{3k+2},\ldots, i_{4k-4}\alpha_{4k-4}).$

For each such $F,$ the number of orderings for the first entry is 2, and the number of orderings for the $\alpha_{h}$ is $[(k-1)!]^{4k}$. Hence, the total number of such $F$ is $2[(k-1)!]^{4k},$ and $$b(F)=(2!)^{4},~c(F)= [(k-1)!]^{4k-8}[(2(k-1))!]^{4},~\pi_{F}(\mathcal{A}_{C_{4,k}})=\big[\cfrac{1}{(k-1)!}\big]^{4k}.$$  Let $D'_{4}$ be the multi-digraph induced by $E(F).$ By the Matrix-Tree Theorem \cite{CvRS}, we have $\tau(D'_{4})=16k^{4k-9}.$ By BEST Theorem \cite{AB}, we have
$\varepsilon(D'_{4})=16k^{4k-9}[(k-2)!]^{4k-8}[(2(k-1)-1)!]^{4}.$
Because each Eulerian cycle has $4k(k-1)$ arcs and multi arcs of $W(F)$ are not labelled, we have
\begin{align*}
\lvert W(F)\rvert&=4k^{4k-8}(k-1)[(k-2)!]^{4k-8}[(2(k-1)-1)!]^{4}.
\end{align*}
For all such $F,$ we know that the total contribution to $\Tr_{4k}(\mathcal{A}_{C_{4,k}})$ is
\begin{align*}
&(k-1)^{n-1}\cfrac{2[(k-1)!]^{4k}(2!)^{4}}{[(k-1)!]^{4k-8}[(2(k-1))!]^{4}}\big[\cfrac{1}{(k-1)!}\big]^{4k}4k^{4k-8}(k-1)
\\&[(k-2)!]^{4k-8}[(2(k-1)-1)!]^{4}
=8k^{4k-8}.
\end{align*}

Therefore, we have \begin{align*}
\Tr_{4k}(\mathcal{A}_{C_{4,k}})=4k^{k-1}(k-1)^{3k-4}+56k^{2k-3}(k-1)^{2k-3}+64k^{3k-5}(k-1)^{k-2}+40k^{4k-8}.
\end{align*}
\hfill$\square$

\noindent\begin{theorem}\label{th:3-5}
Let $C_{4,k}$ be a $k(\geq 4)$-uniform hypercycle with length four. Then
$$\phi_{C_{4,k}}(\lambda)=\lambda^{m_{0}}(\lambda^{k}-1)^{m_{1}}(\lambda^{k}-2)^{m_{2}}(\lambda^{k}-4)^{m_{4}}(\lambda^{k}-\frac{3+\sqrt{5}}{2})^{m'}(\lambda^{k}-\frac{3-\sqrt{5}}{2})^{m'},$$
where $$m_{0}=4(k-1)^{4k-4}-4k^{k-1}(k-1)^{3k-4}+4k^{2k-3}(k-1)^{2k-3}-4k^{3k-5}(k-1)^{k-2}+5k^{4k-8},$$ $$m_{1}=4k^{k-2}(k-1)^{3k-4}-8k^{2k-4}(k-1)^{2k-3}+4k^{3k-6}(k-1)^{k-2}, $$ $$m_{2}=4k^{2k-4}(k-1)^{2k-3}-8k^{3k-6}(k-1)^{k-2}+10k^{4k-9},$$ $$ m_{4}=k^{4k-9}, ~m'=4k^{3k-5}(k-1)^{k-2}-8k^{4k-9}. $$
\end{theorem}

\noindent\textbf{Proof.} By Lemma \ref{le:4-5}, we know the $k$-power of all different eigenvalues of $C_{4,k}$ are $0,~1,~2,~4,\\\frac{3+\sqrt{5}}{2},\frac{3-\sqrt{5}}{2}.$ By Lemma \ref{le:4-4}, we know the trace is a real number. So we know that the multiple of $\frac{3+\sqrt{5}}{2}$ and $\frac{3-\sqrt{5}}{2}$ are equal. Let $m_{0}, m_{1}, m_{2}, m_{4},m'$ be the multiplicities of 0,1,2,4, $\frac{3+\sqrt{5}}{2},$ respectively. By Theorems \ref{th:3-1}-\ref{th:3-4} and $$\Tr_{k}(\mathcal{A}_{C_{4,k}})=k(m_{1}+2m_{2}+4m_{4}+3m'), \Tr_{2k}(\mathcal{A}_{C_{4,k}})=k(m_{1}+4m_{2}+16m_{4}+7m'),$$ $$\Tr_{3k}(\mathcal{A}_{C_{4,k}})=k(m_{1}+8m_{2}+64m_{4}+18m'), ~\Tr_{4k}(\mathcal{A}_{C_{4,k}})=k(m_{1}+16m_{2}+256m_{4}+47m'), $$ $$m_{0}+k(m_{1}+m_{2}+m_{4}+2m')=4(k-1)^{4k-4},$$ we have $$m_{0}=4(k-1)^{4k-4}-4k^{k-1}(k-1)^{3k-4}+4k^{2k-3}(k-1)^{2k-3}-4k^{3k-4}(k-1)^{k-2}+5k^{4k-8},$$ $$m_{1}=4k^{k-2}(k-1)^{3k-4}-8k^{2k-4}(k-1)^{2k-3}+4k^{3k-5}(k-1)^{k-2},$$ $$m_{2}=4k^{2k-4}(k-1)^{2k-3}-8k^{3k-5}(k-1)^{k-2}+10k^{4k-9},$$ $$ m_{4}=k^{4k-9}, ~m'=4k^{3k-6}(k-1)^{k-2}-8k^{4k-9}.$$ \hfill$\square$

\section*{Declaration of competing interest}
The authors declare that they have no conflict of interest.

\section*{Data availability}
The study has no associated data.

\end{document}